\documentclass[12pt,a4paper]{article}

\usepackage{latexsym}
\usepackage{amsfonts}
\usepackage{amssymb}

\newcommand{\be}{\begin{equation}}
\newcommand{\ee}{\end{equation}}
\newcommand{\bea}{\begin{eqnarray}}
\newcommand{\eea}{\end{eqnarray}}
\newcommand{\nn}{\nonumber\\}

\def\G{{\cal G}}                                 %
\def\F{{\cal F}}                                 %
\def\R{{\cal R}}                                 %
\def\D{{\cal D}}                                 %
\def\P{{\cal P}}                                 %
\def\H{{\cal H}}                                 %
\def\I{{\cal I}}                                 %
\def\E{{\cal E}}                                 %
\def\dt{\left.{d\over dt}\right\vert_{t=0}}      %
\def\ad{{\mathrm{ad}}}                           %
\def\Ad{{\mathrm{Ad}}}                           %
\def\half{{\frac{1}{2}}}                         %
\def\quarter{{\frac{1}{4}}}                      %
\begin{document}

\vspace*{0.5cm}
\begin{center}
{\Large \bf On a Poisson-Lie analogue of the classical dynamical 
Yang-Baxter equation for self-dual Lie algebras}
\end{center}

\vspace{1.0cm}

\begin{center}
L. Feh\'er${}^1$ and I. Marshall${}^2$ \\

\bigskip

{\em 
${}^1$Department of Theoretical Physics,
 University of Szeged \\
Tisza Lajos krt 84-86, H-6720 Szeged, Hungary \\
E-mail: lfeher@sol.cc.u-szeged.hu\\

\bigskip

${}^2$
Department of Mathematics, EPFL \\
1015 Lausanne, Switzerland\\
E-mail: ian.marshall@epfl.ch
}
\end{center}

\vspace{1.5cm}

\begin{abstract}
We derive a generalization of the classical dynamical Yang-Baxter 
equation (CDYBE) on a self-dual Lie algebra $\cal G$ by replacing 
the cotangent bundle $T^*G$ in a geometric interpretation
of this equation by its Poisson-Lie (PL) analogue associated with  
a factorizable constant $r$-matrix on $\cal G$.
The resulting PL-CDYBE, with variables in the Lie group $G$
equipped with the Semenov-Tian-Shansky Poisson bracket based on 
the constant $r$-matrix, 
coincides with an equation that appeared in an earlier 
study of PL symmetries in the WZNW model.
In addition to its new group theoretic interpretation,
we present a self-contained analysis of those solutions of the PL-CDYBE  
that were found in the WZNW context and 
characterize them by means of a uniqueness result under a certain 
analyticity assumption.
\end{abstract}




\newpage

\section{Introduction}

The dynamical Yang-Baxter equation plays an 
important role in current research on low dimensional
integrable systems and in related areas of mathematics
(see the reviews \cite{ES,E} and references therein).
In the original form of this equation \cite{Felder} the `dynamical
variable' belongs to the dual of a Cartan subalgebra of a 
simple Lie algebra, $\G$.
Etingof and Varchenko \cite{EV} introduced a generalization
of the classical dynamical Yang-Baxter equation (CDYBE) for which
the variable lies in the dual of an arbitrary subalgebra $\H\subset\G$.
The CDYBE on $\H\subset \G$ ensures the Jacobi identity of
a Poisson bracket (PB) of a certain form on the phase space 
$\H^*\times G\times \H^*$, where $G$ is a connected Lie group
with Lie algebra $\G$, which  thus becomes a 
Poisson-Lie (PL) groupoid \cite{EV}.
The distinguished special case $\H=\G$ 
has interesting applications \cite{AM,BFP}.
The corresponding PB on $\G^* \times G\times \G^*$ can be viewed 
as a modification of the canonical PB of the cotangent bundle
$T^*G$ expressed in terms of redundant variables.

The idea developed in the present Letter is to replace $T^*G$  
in the construction of \cite{EV} by its PL analogue,
the Heisenberg double \cite{STS,AF}, associated with
a factorizable constant $r$-matrix.
This will lead to a natural PL analogue of the CDYBE on $\G$.
In fact, we shall obtain equation (\ref{21}), which we call the `PL-CDYBE',  
as the guarantee of the
Jacobi identity of a PB 
on a PL groupoid of the form $\check G\times G \times \check G$,
where $\check G$ is a neighbourhood of $e\in G$ diffeomorphic to
a corresponding domain in the dual $G^*$ of the PL group $G$
equipped with the quadratic $r$-bracket.

The PL-CDYBE turns out to 
coincide with a special case of the generalized CDYBE 
that governs the chiral WZNW PBs for generic monodromy \cite{BFP}.
This special case corresponds to PL symmetries in the chiral
WZNW model, and in this context some 
remarkable solutions of (\ref{21}) have been found in \cite{BFP}, too, 
although the fact that they are solutions was proved in quite
an indirect manner only.
We here give a self-contained presentation of these solutions,
which has the advantage of showing also their uniqueness 
under a certain analyticity assumption.
It will be clear that these solutions are the natural PL analogues 
of the canonical (Alekseev-Meinrenken) solution of the
CDYBE on $\G$ found in \cite{EV,AM,BFP}.

The present Letter may serve as a starting point towards 
deriving generalizations of the CDYBE on $\H\subset \G$ 
for which the dynamical 
variable belongs to the dual $H^*$ of a PL subgroup $H$ of a PL group $G$.
For this, one should extend the observation \cite{FGP} that 
the CDYBE on $\H\subset \G$ can be 
obtained by applying Dirac reduction to the
PL groupoid that encodes the CDYBE on $\G$.
We hope to return to this question in the future.

As for the organization of the rest of the Letter,
the derivation of the PL-CDYBE is contained in Section 2,
its analytic solutions are described in Section 3, and are
further discussed in Section 4.
Our results are summarized by Proposition 1 in Section 2 and
by Proposition 2 in Section 3.

\section{PL generalization of the CDYBE on $\G$}

Consider a (real or complex) self-dual 
Lie algebra\footnote{A Lie algebra that admits an
invariant scalar product, i.e. an invariant nondegenerate symmetric bilinear form,
is called self-dual.  
A review of self-dual  Lie algebras can be found in \cite{Figu}.}  
$\G$ with 
a fixed invariant scalar product $\langle\ ,\ \rangle$ and  
a corresponding connected Lie group $G$.
The scalar product is used to identify
the dual space $\G^*$ with $\G$.
We below recall the derivation of the CDYBE on $\G$ and then generalize
it to the PL case in correspondence with any antisymmetric 
constant solution $r\in \G\otimes \G$ of the modified CYBE:
\be
[r_{12},r_{13}]+\hbox{cycl. perm.}=-\frac{1}{4} f,
\label{1}\ee
where $f= f_{ab}^c T^a \otimes T^b\otimes T_c$ with 
dual bases $\{T_a\}$ and $\{T^b\}$ of $\G$  satisfying 
$[T_a,T_b]=f_{ab}^c T_c$ and $\langle T_a, T^b\rangle=\delta_a^b$.
Summation over coinciding indices is understood and the notations
$C:= T_a \otimes T^a$, $r^\pm := r \pm \frac{1}{2} C$ are used 
throughout the Letter.

\subsection{Recall of the CDYBE on $\G$}

Let us introduce the phase space 
\be
{\cal P}:= \G \times G \times \G = 
\{ (\omega^L, g, \omega^R)\,\vert\, \omega^{L,R}\in \G,\, g\in G \}
\label{2}\ee
and equip it with the Poisson structure $\{\ ,\ \}_\P^0$ defined as follows:
\bea
&&\{ g_1, g_2\}_\P^0= 0\nn
&& \{ \omega^L_1, \omega^L_2\}_\P^0= \half [C, \omega^L_2 -
\omega^L_1] \nn
&& \{ \omega^L_1,g_2\}_\P^0 = C g_2\nn
&& \{ \omega^R_1, \omega^R_2\}_\P^0=-\half
[C, \omega^R_2 -\omega^R_1]\nn
&& \{ \omega^R_1, g_2\}_\P^0 = g_2 C\nn
&&\{ \omega^R_1, \omega^L_2\}_\P^0=0.
\label{3}\eea
In this description of the PBs of 
the fundamental variables $g$, $\omega^{L,R}$ we employ
the standard tensorial notation \cite{FadTak} that implicitly refers to
an arbitrary matrix representation of $G$.
The formulae  
\be
\phi: G \times \P \rightarrow \P,
\quad
\phi: (q, (\omega^L, g, \omega^R))\mapsto 
(q \omega^L q^{-1}, qg, \omega^R)
\label{4}\ee
\be
\psi:\P \times G \rightarrow \P,
\quad
\psi: ((\omega^L, g, \omega^R), q)\mapsto 
(\omega^L , gq, q^{-1}\omega^R q)
\label{5}\ee
define left- and right  actions of the group $G$ on $\P$.
These are Poisson actions in the usual sense,
with respect to the zero PB on $G=\{q\}$, 
and the corresponding momenta that generate them are directly given
by $\omega^L$ and $\omega^R$, respectively. 
The constraint 
\be
\omega^R = g^{-1} \omega^L g
\label{6}\ee
defines a Poisson submanifold of $(\P, \{\ ,\ \}_\P^0)$, which can 
be identified with the cotangent bundle $T^*G$  with its canonical PB.
Conversely,  $(\P, \{\ ,\ \}_\P^0)$ is obtained from the 
cotangent bundle by `forgetting' this familiar relation 
between the left- and right momenta.

Now consider a (smooth or holomorphic) map 
$\R: \check \G \rightarrow \G\wedge \G$, where $\check \G \subset \G$
is an open submanifold stable under the adjoint action of $G$.
By using $\R$, let us try to define a PB $\{\ ,\ \}_{\check \P}$
on the manifold $\check \P:= \check \G \times G \times \check \G\subset \P$
by the following formula:
\bea
&&\{ g_1, g_2\}_{\check \P}= \R(\omega^L) g_1 g_2 - g_1 g_2 \R(\omega^R)\nn
&& \{ \omega^L_1, \omega^L_2\}_{\check \P}= 
\half [C, \omega^L_2 -\omega^L_1] \nn
&& \{ \omega^L_1,g_2\}_{\check \P} = C g_2\nn
&& \{ \omega^R_1, \omega^R_2\}_{\check \P}=
-\half [C, \omega^R_2 -\omega^R_1]\nn
&& \{ \omega^R_1, g_2\}_{\check \P} = g_2 C\nn
&&\{ \omega^R_1, \omega^L_2\}_{\check \P}=0.
\label{7}\eea
This formula differs from (\ref{3}) only in the first line that
contains $\R$.

It is easy to check that the Jacobi identities of 
$\{\ ,\ \}_{\check \P}$ require $\R$ to be equivariant,
\be
\R(q \omega q^{-1})= (q\otimes q) \R(\omega) (q^{-1}\otimes q^{-1})
\qquad
\forall q\in G,\, \omega=\omega^aT_a \in \check \G,
\label{8}\ee  
and to satisfy the equation 
\be
[\R_{12}, \R_{13}] + T^a_3 \frac{\partial}{\partial \omega^a} \R_{12} 
+\hbox{cycl. perm.} = {\I} \qquad\hbox{on}\qquad \check \G,
\label{9}\ee
where $\I$ is some $G$-invariant constant element 
of $\G\wedge \G\wedge \G$. 
Equation (\ref{9}) is the (modified) CDYBE on $\G$.
(Strictly speaking, the term `modified' should be 
used if $\I\neq 0$.)

The above interpretation of the CDYBE on $\G$ is taken from \cite{EV}.
Due to the equivariance property (\ref{8}),
the actions of $G$ given by $\phi$ (\ref{4}) and 
$\psi$ (\ref{5}) yield Poisson actions
on $(\check \P, \{\ ,\ \}_{\check \P})$ that are generated by the
momenta $\omega^{L,R}$ in the same way as in the $\R=0$ case.
It is also explained in \cite{EV} that $(\check \P, \{\ ,\ \}_{\check \P})$ 
has the structure of a PL groupoid. 

Let $\F$ denote the holomorphic complex function given by
\be
\F(z)=  \coth z - z^{-1},
\qquad
\F(0)=0.
\label{can1}\ee
By using the identification $\G\otimes \G \simeq {\mathrm{End}}({\G})$
defined by the scalar product on $\G$, one 
obtains \cite{EV,AM,BFP} a solution 
of (\ref{9}) by setting 
\be
\R(\omega):= \F(\ad_\omega),
\qquad\qquad
\I=-f 
\label{can2}\ee
with $f$  in (\ref{1}). 
In \cite{FP} this $r$-matrix is called `canonical' since it is actually 
the unique solution of (\ref{9}) under assuming 
(\ref{can2}) as an ansatz with {\em some} holomorphic odd function 
regular at $0$. 
Generalizations of the canonical  
$r$-matrix are described in Section 3.

\subsection{Derivation of the PL-CDYBE}

A natural PL analogue, $(P, \{\ ,\ \}_P^0)$, of the Poisson 
manifold $(\P, \{\ ,\ \}_\P^0)$ is provided  by
\be
P:= G \times G \times G = 
\{ (\Omega^L, g, \Omega^R)\,\vert\, \Omega^{L,R}\in G,\, g\in G \}
\label{10}\ee
with the PB $\{\ ,\ \}_P^0$ defined by 
\bea
&&\{ g_1, g_2\}_P^0= r g_1 g_2 - g_1 g_2 r \nn
&& \{ \Omega^L_1, \Omega^L_2\}_P^0= r \Omega^L_1 \Omega^L_2 +
\Omega^L_1 \Omega^L_2 r  - 
\Omega^L_1 r^- \Omega^L_2 - \Omega^L_2 r^+ \Omega^L_1\nn
&& \{ \Omega^L_1,g_2\}_P^0 = ( r^+ \Omega^L_1 - \Omega^L_1 r^-) g_2\nn
&& \{ \Omega^R_1, \Omega^R_2\}_P^0=
 - \left(r \Omega^R_1 \Omega^R_2 + \Omega^R_1 \Omega^R_2 r 
 -\Omega^R_1 r^- \Omega^R_2 - 
\Omega^R_2 r^+ \Omega^R_1\right)\nn
&& \{ \Omega^R_1, g_2\}_P^0 = g_2 ( r^+ \Omega^R_1 - 
\Omega^R_1 r^-)\nn
&&\{ \Omega^R_1, \Omega^L_2\}_P^0=0.
\label{11}\eea
The Poisson submanifold $(T^* G)_r \subset P$ specified by the constraint
\be
\Omega^R = g^{-1} \Omega^L g
\label{12}\ee
is the PL analogue of the cotangent bundle considered in \cite{STS, AF}.
The formulae
\be
\Phi: G \times P \rightarrow P,
\quad
\Phi: (q, (\Omega^L, g, \Omega^R))\mapsto 
(q \Omega^L q^{-1}, qg, \Omega^R)
\label{13}\ee
\be
\Psi:P \times G \rightarrow P,
\quad
\Psi: ((\Omega^L, g, \Omega^R), q)\mapsto 
(\Omega^L , gq, q^{-1}\Omega^R q)
\label{14}\ee
define left- and right  actions of the group $G$ on $P$.
These are PL actions if the group 
$G=\{q\}$ is endowed with the Sklyanin PB
\be
\{ q_1, q_2\} = r q_1 q_2 - q_1 q_2 r, 
\label{15}\ee
and  then $\Omega^{L,R}$ yield group valued momenta 
generating  these actions.
It is well known \cite{STS} that the PB on $G$ 
appearing in the second line of
(\ref{11})  becomes the natural PB on the PL group $G^*$ dual to
$G$ if we locally identify $G^*$ with $G$ in the standard manner.
This PB on $G$ is sometimes referred to as the Semenov-Tian-Shansky PB.
The Poisson space $(\P, \{\ ,\ \}_\P^0)$ is 
a linearization of
$(P,\{\ ,\ \}_P^0)$ in the same sense in which $T^*G$ is the linearization
of its PL analogue (see e.g. \cite{AF}).  
This is further discussed in Section 4.

Now it should be clear how the above set of analogies
can be extended to obtain a natural PL generalization of the CDYBE on $\G$.
As input data, we consider a (smooth or holomorphic) map
$R: \check G \rightarrow \G\wedge \G$, where $\check G\subset G$
is an open submanifold stable under conjugation by any $q\in G$.
We then require $R$ to define a PB,  $\{\ ,\ \}_{\check P}$, 
on the manifold 
$\check P:= \check G \times G \times \check G$ by means of following 
modification of (\ref{11}):
\bea
&&\{ g_1, g_2\}_{\check P}= (r+R(\Omega^L)) g_1 g_2 
- g_1 g_2 (r+R(\Omega^R)) \nn
&& \{ \Omega^L_1, \Omega^L_2\}_{\check P}= r \Omega^L_1 \Omega^L_2 +
\Omega^L_1 \Omega^L_2 r  - 
\Omega^L_1 r^- \Omega^L_2 - \Omega^L_2 r^+ \Omega^L_1\nn
&& \{ \Omega^L_1,g_2\}_{\check P} = ( r^+ \Omega^L_1 - 
\Omega^L_1 r^-) g_2\nn
&& \{ \Omega^R_1, \Omega^R_2\}_{\check P}=
 - \left(r \Omega^R_1 \Omega^R_2 + \Omega^R_1 \Omega^R_2 r  
-\Omega^R_1 r^- \Omega^R_2 - 
\Omega^R_2 r^+ \Omega^R_1\right)\nn
&& \{ \Omega^R_1, g_2\}_{\check P} = g_2 ( r^+ \Omega^R_1 - 
\Omega^R_1 r^-)\nn
&&\{ \Omega^R_1, \Omega^L_2\}_{\check P}=0.
\label{16}\eea
We are interested in such Poisson structures for which the restrictions of 
$\Phi$ (\ref{13}) and $\Psi$ (\ref{14}) to $\check P$ yield PL actions
of the group $G$ equipped with the Sklyanin PB (\ref{15}).
As is easy to check, the ansatz (\ref{16}) enjoys this PL symmetry 
if and only if   $R$ is equivariant,
\be
R(q \Omega q^{-1})= (q\otimes q) R(\Omega) (q^{-1}\otimes q^{-1})
\qquad
\forall q\in G,\, \Omega\in \check G.
\label{17}\ee  
Therefore we assume (\ref{17}) to hold.
The equivariance of $R$ also guarantees the Jacobi identities
\be
\{\{ g_1, g_2\}_{\check P}, \Omega^L_3\}_{\check P} 
+ \hbox{cycl. perm}=0,
\qquad
\{\{ g_1, g_2\}_{\check P}, \Omega^R_3\}_{\check P} + \hbox{cycl. perm}=0.
\label{18}\ee 
To explain how this works, 
let $\lambda^a$ and $\rho^a$ denote the vector fields on $G$ 
that operate on the group element by left- and by right multiplication 
by $T^a\in \G$, respectively. In other words, $\lambda^a$ and $\rho^a$
are the right- and left-invariant vector fields associated
with $T^a\in T_eG$, respectively. Introduce also 
\be
\D^a_\pm := (\rho^a \pm \lambda^a)
\qquad\hbox{and}\qquad
\D_a^\pm := (\rho_a \pm \lambda_a).
\label{*}\ee
In fact, (\ref{16}) implies the relation  
\bea
&&\{\{ g_1, g_2\}_{\check P}, \Omega^L_3\}_{\check P} +\hbox{cycl. perm.}=
\left( [R_{12}(\Omega^L), r^-_{13}+ r^-_{23}] 
-r^-_{ab} \D_-^a R_{12}(\Omega^L) T^b_3 \right)\Omega^L_3 g_1 g_2\nn
&&\qquad\qquad - \Omega^L_3 \left( [R_{12}(\Omega^L), r^+_{13}+ r^+_{23}] 
-r^+_{ab} \D_-^a R_{12}(\Omega^L) T^b_3 \right)g_1 g_2.
\label{19}\eea
The right hand side vanishes as a consequence of the infinitesimal 
version of (\ref{17}), 
\be
\D^a_- R_{12}(\Omega) = [R_{12}(\Omega), T^a_1 + T^a_2] 
\qquad
\forall \Omega\in \check G,\, T^a\in \G.
\label{20}\ee
The other Jacobi identity in (\ref{18}) holds in a similar manner.
The Jacobi identities involving 
$\{\{ g_1, \Omega^L_2\}_{\check P}, \Omega^L_3\}_{\check P}$ and
$\{\{ g_1, \Omega^R_2\}_{\check P}, \Omega^R_3\}_{\check P}$ 
are automatically satisfied since 
they do not contain $R$. 
Note that $\Omega^{L,R}$ serve as momentum maps for
$(\check P, \{\ ,\ \}_{\check P})$ in the same way 
as they do for $(P, \{\ ,\ \}_P^0)$.
The only further requirement on $R$ imposed by the Jacobi identities of 
$\{\ ,\ \}_{\check P}$ is the PL-CDYBE described in following proposition.

\bigskip\noindent
{\bf Proposition 1.}
{\em Let $R: \check G\rightarrow \G\wedge \G$ be a 
(smooth or holomorphic) $G$-equivariant map (\ref{17}). 
Formula (\ref{16}) generates a PB on the (smooth or holomorphic)
functions on $\check P$ if and only if  $R$ 
satisfies the PL-CDYBE given by}
\be
[R_{12}, R_{13}] + \half T^a_3 \D_a^+ R_{12} 
+\hbox{cycl. perm.} = {\I} \qquad\hbox{on}\qquad \check G,
\label{21}\ee
{\em 
where $\cal I$ is an arbitrary $G$-invariant constant element 
of $\G\wedge \G\wedge \G$.} 
 
\bigskip\noindent
{\em Proof.}
It is clear that (\ref{16}) defines a PB if and only if
the Jacobi identity involving 
$\{\{g_1, g_2\}_{\check P}, g_3 \}_{\check P}$ is
satisfied. 
A straightforward calculation gives
\be
\{ \{ g_1, g_2\}_{\check P}, g_3\}_{\check P} 
+ \hbox{cycl. perm.}= \E(\Omega^L)(g_1 g_2 g_3)
-(g_1 g_2 g_3)\E(\Omega^R),
\label{22}\ee
where the expression ${\cal E}$ reads as
\be
\E(\Omega)=\left([r_{12}+R_{12},r_{13}+R_{13}] + 
T^a_3 (r_{ab}{\cal D}^b_- +\half {\cal D}_a^+)
 R_{12} +\hbox{cycl. perm.}\right)(\Omega).
\label{23}\ee
By (\ref{22}),  the Jacobi identity requires 
$\E(\Omega)$ to be a $G$-invariant constant.
We notice that the equivariance (\ref{20}) of $R$ implies 
the relation   
\be
[r_{12}+ R_{12},r_{13}+R_{13}]+\hbox{cycl. perm.}
=[r_{12},r_{13}] + [R_{12},R_{13}] - 
T_3^a r_{ab} {\cal D}_-^b R_{12}
+\hbox{cycl. perm.}
\label{24}\ee 
By inserting  this and (\ref{1}) into (\ref{23}) we immediately
obtain the PL-CDYBE (\ref{21}). {\em Q.E.D.}

\medskip\noindent
{\em Remark 1.} 
It is remarkable that the PL-CDYBE (\ref{21}) does not
contain any explicit reference to the background constant $r$-matrix  
used in (\ref{16}). 
This equation first appeared in  \cite{BFP} as a condition 
on the `exchange $r$-matrices' associated with
PL symmetries in the chiral WZNW model.
In this application the role of the group valued variable  
is played by the monodromy matrix of the chiral WZNW field, but the
PBs of the monodromy matrix  are different from the Semenov-Tian-Shansky PBs 
that appear in the second and fourth lines of (\ref{16}).

\medskip\noindent {\em Remark 2.}  
The manifold $\check P= \check G \times G \times \check G$
is a groupoid \cite{Mackan} with the partial multiplication  
\be
(\bar \Omega^R, \bar g, \bar \Omega^L) (\Omega^R, g, \Omega^L)
= (\hat \Omega^R, \hat g, \hat \Omega^L) 
\label{G1}\ee
defined by the constraints
\be
\bar \Omega^L= \Omega^R,
\quad
\hat \Omega^R= \bar \Omega^R,
\quad
\hat \Omega^L= \Omega^L,
\quad
\hat g= \bar g g.
\label{G2}\ee
These are  first class constraints on $\check P \times \check P 
\times \check P^-$ if the first two factors are equipped with the PB 
$\{\ ,\ \}_{\check P}$ and $\check P^-$, 
the set of the hatted triples in (\ref{G1}),
is equipped with the opposite PB.
This means that the phase space $(\check P, \{\ ,\ \}_{\check P})$ 
is a PL groupoid \cite{Weinstein}.
As far as the PBs are concerned, 
this PL groupoid appears to be different from the groupoids
associated with the  
WZNW exchange $r$-matrices in \cite{BFP} even in the special case when
those exchange $r$-matrices satisfy the PL-CDYBE (\ref{21}). 
It would be very interesting to clarify the relationship 
(which is perhaps an equivalence by some change of variables) between 
the PL groupoids $(\check P, \{\ ,\ \}_{\check P})$ and those 
constructed in \cite{BFP}.

\section{A family of solutions of the PL-CDYBE}

To simplify some arguments, 
in the main body of this section we assume 
that $\G$ is a complex simple Lie algebra and describe certain solutions 
of (\ref{21}) in this case.
These dynamical $r$-matrices were originally found in the context of 
the WZNW model \cite{BFP}.
Their presentation below is considerably simpler than the one in \cite{BFP} 
and our Proposition 2 includes a new uniqueness result as well.  
At the end of the section, it will be remarked that these solutions 
of (\ref{21}) are available for any self-dual Lie algebra.

First note that for a simple Lie algebra 
the invariant $\I$ in (\ref{21}) must have the form 
\be
\I = \mu f,
\label{25}\ee
where $f$ is given in (\ref{1}) and $\mu$ is some constant.
If (\ref{25}) holds and we 
identify $\G\otimes \G$ with ${\mathrm{End}}(\G)$ by
the  scalar product
on $\G$, then (\ref{21}) can be equivalently rewritten as the requirement 
\be
\langle [ R(\Omega)X,  R(\Omega)Y] 
-\half {\cal D}_X^+ R(\Omega) Y, Z\rangle 
+\hbox{cycl. perm.} = \mu \langle [X, Y],Z\rangle
\quad
\forall X,Y,Z\in\G,
\label{26}\ee 
where the cyclic permutations act on $X,Y,Z$ and we have
\be 
{\cal D}_X^+ R(\Omega)= X^a {\cal D}_a^+ R(\Omega)=
\dt R(e^{t X} \Omega e^{tX})
\quad
\hbox{for} \quad R: \check G\rightarrow {\mathrm{End}}(\G).
\label{27}\ee 
We  suppose the domain $\check G\subset G$ to be diffeomorphic
to a neighbourhood of zero $\check \G\subset \G$ by means of the
exponential map, and thus parametrize $\Omega \in \check G$ 
according to 
\be
\Omega = e^\omega 
\qquad\hbox{with}\quad \omega \in \check \G.
\label{28}\ee
We further write 
\be
R(e^\omega) =\tilde R(\omega) 
\quad\hbox{for}\quad
\omega\in \check \G.
\label{29}\ee
Define the holomorphic complex function $h$ by 
\be
h(z)=\half z \coth (\half z),
\qquad h(0)=1.
\label{30}\ee
In fact, the required derivatives of $R$ can be expressed in terms
of the variable $\omega$ as 
\be
\half {\cal D}_X^+ R(e^\omega) = 
\dt \tilde R(\omega + h (\ad_\omega) X t ).
\label{31}\ee
This is a consequence of standard identities (e.g. \cite{SW}) 
expressing the left- and right-invariant vector fields on $G$
in the exponential parametrization.

Now let $F$ be a holomorphic complex function which 
is regular in a neighbourhood of $z=0$ and is odd, $F(-z)=-F(z)$.
Then consider the ansatz 
\be
R(e^\omega) = F(\ad_\omega),
\qquad
\omega\in \check \G.
\label{32}\ee
Note that $F(\ad_\omega)$ is given by means of   
the Taylor series of $F(z)$ around $z=0$
if $\omega$ is near enough to zero.  
This ansatz automatically guarantees the equivariance and
the antisymmetry of the dynamical $r$-matrix $R(e^\omega)$.

\bigskip\noindent {\bf Proposition 2.}
{\em Let $\G$ be a complex simple Lie algebra.
The ansatz (\ref{32}) provides a solution of the PL-CDYBE (\ref{26}),
which is equivalent to (\ref{21}), on a 
domain $\check G=\exp(\check \G)$
if and only if the holomorphic odd function $F(z)$ is given by
\be
F_\nu(z)=\nu\coth(\nu z) - \half\coth(\half z),
\label{33}\ee
where $\nu$ is an arbitrary constant related to the 
constant $\mu$ in (\ref{25}) by
$\mu=\quarter-\nu^2$.}

\bigskip\noindent {\em Proof.}
We start by noting that 
equation (\ref{26}) is satisfied on $\check \G$ if and only if
it is satisfied on the dense, open submanifold of $\check \G$ 
consisting of regular semisimple elements. 
Therefore we may subsequently fix 
$\omega$ to be an arbitrarily chosen regular
semisimple element of $\check \G$.
For such an $\omega$ it is convenient to 
choose a basis of $\G$ spanned by root vectors 
$E_\alpha$ ($\alpha\in \Delta$), where $\Delta$ is the set of roots
with respect to the Cartan subalgebra $\H$ that contains $\omega$,
and a basis $H_i$ of $\H$.
The ansatz (\ref{32}) gives directly that
\be
R(e^\omega) H_i=0,
\qquad
R(e^\omega) E_\alpha = F(\alpha(\omega)) E_\alpha.
\label{35}\ee
The derivatives appearing in (\ref{26}) can be computed 
 with the aid of (\ref{31}) and the equivariance of $R$,
$R( q \Omega q^{-1}) = \Ad_q \circ R(\Omega) \circ \Ad_{q^{-1}}$
for $\Ad_q\in {\mathrm{End}}(\G)$ $\forall q\in G$.
We find 
\bea 
&&  \D^+_{E_\alpha} R(e^\omega)=2 \theta(\alpha(\omega)) 
[ R(e^\omega), \ad_{E_\alpha}]\qquad \forall \alpha\in \Delta,
\nonumber\\
&&  \D^+_{H_i} R(e^\omega)=2 F'(\ad_\omega) \circ \ad_{H_i}
\qquad \forall H_i\in \H,
\eea
where we introduced the notation 
\be
\theta(z) := z^{-1} h(z)= \half \coth(\half z).
\label{40}\ee
This can be spelled out as 
\bea
&&  \D^+_{E_\alpha}R(e^\omega)E_\beta 
        =2\theta(\alpha(\omega))
 \Bigl(F(\alpha(\omega)+\beta(\omega)) 
      - F(\beta(\omega))\Bigr) [E_\alpha,E_\beta] \nonumber\\
&&  \D^+_{E_\alpha}R(e^\omega)H_i =
      - 2\alpha(H_i) \theta(\alpha(\omega)) 
F(\alpha(\omega)) E_\alpha \nonumber\\
&&  \D^+_{H_i}R(e^\omega)E_\alpha =
       2\alpha(H_i)F'(\alpha(\omega)) E_\alpha \nonumber\\
&&  \D^+_{H_i} R(e^\omega) H_j =0,
\label{37}\eea
for any $\alpha, \beta \in \Delta$.
The first line of (\ref{37}) simplifies if $\beta=-\alpha$, 
since  $F(0)=0$ by the oddness of $F$.
The nontrivial conditions represented by 
(\ref{26}) arise in the cases
\be 
X=E_\alpha, Y= H_i, Z=E_{-\alpha}
\quad \hbox{and}\quad
X=E_\alpha, Y=E_\beta, Z=E_{-\alpha-\beta},
\label{34}\ee
where $H_i\in \H$ and $\alpha, \beta, (\alpha+\beta)\in \Delta$.
Indeed, 
by evaluating (\ref{26}) for the 
independent choices in (\ref{34}) we obtain that 
(\ref{26}) under the ansatz (\ref{32}) is equivalent to the
following conditions on the holomorphic function $F$:
\be
F'(z)+2 \theta(z) F(z)  + F^2(z) +\mu=0,
\qquad F(0)=0,
\label{38}\ee
and 
\bea
&& F(z)F(w) - F(z+w)( F(z) + F(w) ) -
   \theta(z) ( F(z+w) - F(w) ) \nonumber\\
&& \qquad\quad  - \theta(w) ( F(z+w) - F(z) ) -
   \theta(z+w) ( F(z) + F(w) ) -\mu=0.
\label{39}\eea
The differential equation (\ref{38}) arises from the first choice
in (\ref{34}), while 
the functional equation (\ref{39}) arises from the
second choice.
Thus (\ref{39}) must hold for nonzero 
$z$, $w$, $(z+w)$ in a neighbourhood of $0$.

Now the differential equation (\ref{38}) 
is not difficult to solve. 
Let $\chi(z)=F(z)+\theta(z)$. 
Then the identity $\theta^2(z)+\theta'(z)={1\over 4}$
and setting $\mu=\frac{1}{4}-\nu^2$ gives
\be
\chi'(z)+\chi^2(z)=\nu^2,
\qquad\qquad
\lim_{z\rightarrow 0}z\chi(z)=1.
\label{41}\ee
The unique solution of this is 
\be
\chi(z)=\nu\coth(\nu z),
\label{42}\ee
which yields  
\be
F(z)=\nu\coth(\nu z) - \half \coth(\half z),
\qquad
\mu=\quarter-\nu^2.
\label{43}\ee
So far we have shown that 
for (\ref{32}) to satisfy (\ref{26}) it is necessary 
for $F$ to be given by
(\ref{43}) with some $\nu$. 
To see that it is also sufficient we must verify the
functional equation (\ref{39}). 
Using the identity $\coth x\coth y-\coth(x+y)(\coth x+\coth y)=-1$,
it is easy to check that
for $F$ given by (\ref{43}), (\ref{39}) does indeed hold.
{\em Q.E.D.}

\medskip\noindent {\em Remark 3.}
Let us comment on the domain of definition of the 
$r$-matrices provided by Proposition 2.
The map $R: \check G: \rightarrow {\mathrm{End}}(\G)$ 
is defined in (\ref{32}) by a power series in $\ad_{\log \Omega}$ if 
$\Omega$ is near to $e\in G$. 
For generic $\nu$, the domain of this $r$-matrix can be extended
naturally to contain all $\Omega$ that has a unique logarithm and the
eigenvalues of $\ad_{\log \Omega}$ do not intersect the poles
of the holomorphic function $F_\nu$.
If $\nu$ is a half-integer, then our $r$-matrix can be 
expressed directly in terms of $\Omega$.
For example, we obtain
\be
R(\Omega) = \half \left(\Ad_\Omega -1\right)(\Ad_\Omega +1)^{-1} 
\qquad 
\hbox{for}\quad \nu=1,
\label{nu1}\ee
since $F_1(z) = \half (e^z -1)(e^z +1)^{-1}$.
The maximal domain of definition of $R$ in (\ref{nu1})  
contains all $\Omega \in G$ for which $-1$ is not an eigenvalue
of $\Ad_\Omega$.
It is easy to verify directly that 
(\ref{nu1}) solves the PL-CDYBE (\ref{26}) 
with $\mu=-\frac{3}{4}$ for any
Lie group $G$ with a self-dual Lie algebra.

\medskip\noindent {\em Remark 4.}
In fact, the statement of Proposition 2 remains valid 
if we replace the complex simple Lie algebra $\G$ with any 
(complex or real) self-dual Lie algebra.
This can be proved by the method used in 
\cite{FP} to analyse  the CDYBE (\ref{9}) on an 
arbitrary self-dual Lie algebra.
The arguments contained in \cite{BFP} are also valid in this generality.

\section{Discussion}

In this Letter we have shown that equation (\ref{21}) is a natural
PL analogue of the CDYBE (\ref{9}) on $\G$ and established
a uniqueness result concerning the family of 
solutions of it  given by (\ref{32}) with (\ref{33}). 
The PL-CDYBE together with these solutions appeared  
earlier in the concrete context of PL symmetries on
the chiral WZNW phase space \cite{BFP}.
The  present work provides a new,  
purely group theoretic interpretation of this equation.  

Since $T^*G$ can be viewed as a `scaling limit' of $(T^*G)_r$ \cite{AF}, 
it should be possible to view also the canonical solution (\ref{can2}) 
of the CDYBE on $\G$ as a limiting case of the
$r$-matrices given in Proposition 2. 
We conclude by outlining how this comes about.

Let $\{\ ,\ \}_{\check P, \nu}$ denote 
the PB (\ref{16}) on $\check P$ 
defined by $R(\Omega):= F_\nu(\ad_{\log \Omega})$ with (\ref{33}).
For an arbitrary constant $\gamma$, consider the map 
$\pi_\gamma: \check \P \rightarrow \check P$ given by
the scaled exponential parametrization 
\be
\pi_\gamma: (\omega^L, g, \omega^R)\mapsto 
 (e^{\gamma \omega^L}, g, e^{\gamma \omega^R}).
\label{D1}\ee
By choosing the domains $\check \G$ and $\check G$ so that this map
is one-to-one\footnote{Of course, $\omega^{L}$ and $\omega^R$ 
are also restricted by the 
condition that $F_\nu(\gamma \ad_\omega)$ must 
be well-defined for $\omega=\omega^L, \omega^R$,
but here we suppress this dependence of the domain $\check \P$ on 
the parameters $\gamma$ and $\nu$ for brevity.   
A maximal domain arises naturally
in the limit $\gamma \rightarrow 0$, which is used in (\ref{D3}).},  
the PB $\gamma \{\ ,\ \}_{\check P, \nu}$ on $\check P$
gives rise to a PB on $\check \P$, which we denote by
$\{\ ,\ \}_{\check \P}^{\nu,\gamma}$.
We then perform the limit $\gamma \rightarrow 0$, but in such a way
that at the same time we let $\nu$  to be a function $\nu_\gamma$ 
for which
\be
\lim_{\gamma \rightarrow 0} (\gamma \nu_\gamma) = \tau
\label{D2}\ee
with some constant $\tau$. 
For functions $\Psi_1, \Psi_2$, the formula
\be
\{ \Psi_1, \Psi_2\}_{\check \P}^\tau(\omega^L,g,\omega^R)
:= \lim_{\gamma \rightarrow 0} 
\{\Psi_1,\Psi_2 \}_{\check \P}^{\nu_\gamma,\gamma}(\omega^L,g,\omega^R)
\label{D3}\ee
is easily seen to define a PB on $\check \P$. 
In fact, by inspecting the PBs of the basic variables  
$g, \omega^{L},\omega^{R}$,  we find that the PB $\{\ ,\ \}_{\check \P}^\tau$
has the form of the ansatz (\ref{7}) with 
\be
\R(\omega)= \F_\tau(\ad_\omega)
\qquad\hbox{where}\qquad
\F_\tau(z):= \tau \coth(\tau z) - z^{-1}.
\label{D4}\ee
Since $\F_\tau=0$ for $\tau=0$,  $\{\ ,\ \}_{\check \P}^\tau$
reproduces the original PB (\ref{3}) on $\check \P$ for $\tau=0$.  
For $\tau=1$ (\ref{D4}) becomes the canonical $r$-matrix in (\ref{can2}),
and for an arbitrary $\tau\neq 0$ it yields a solution of the
CDYBE (\ref{9}) with right hand side given by $\I = - \tau^2 f$.

Thus the canonical $r$-matrix (\ref{can2}) 
is indeed a scaling limit
of the $r$-matrices given by Proposition 2.
This $r$-matrix serves as a source of many 
solutions of the CDYBE on subalgebras $\H\subset \G$
by using Dirac reduction and by taking various limits \cite{EV,FGP}. 
Its generalizations given by Proposition 2 may perhaps play an
analogous role in the PL context, but this leads 
to questions outside the scope of the present Letter.
Another open problem that could be interesting to study
is to quantize these $r$-matrices. 
The related problem of quantizing solutions
of the CDYBE on a non-Abelian base, such as on $\G$,
has been investigated recently in \cite{Xu}.

\medskip
\bigskip
\noindent
{\em Note added in proof upon publication in Lett.~Math.~Phys.:}  In the meantime 
we have found the change of variables alluded to at the end of Section 2.  
As will be described elsewhere, it operates by mapping the group-valued variable  
$\Omega$ to the WZNW monodromy matrix $M:= \exp(2\nu \log \Omega)$.

\section*{Acknowledgments}
We are grateful to A. G\'abor and J. Balog for discussions and for comments
on the manuscript.
L.F. was supported in part by the Hungarian 
Scientific Research Fund (OTKA) under T034170, T030099, T029802 and M036804.

\end{document}